\documentclass{amsart}
\usepackage{amsmath}
\usepackage{amstext}
\usepackage{amssymb}
\usepackage{amsthm}
\swapnumbers
\theoremstyle{plain}
\newtheorem{thm}{Theorem}[section]
\newtheorem{lem}[thm]{Lemma}

\newtheorem{cor}[thm]{Corollary}

\theoremstyle{definition}
\newtheorem{defi}[thm]{Definition}
\newtheorem{defs}[thm]{Definitions}

\newtheorem{ntn}[thm]{Notation}
\newtheorem{rmd}[thm]{Reminder}

\theoremstyle{remark}

\newtheorem{rmk}[thm]{Remark}

 \DeclareMathOperator{\height}{ht}

 \DeclareMathOperator{\ass}{ass}

 \DeclareMathOperator{\ann}{ann}

\def\Z{\mathbb Z}
\def\N{\mathbb N}

\def\fa{{\mathfrak{a}}}

\def\fb{{\mathfrak{b}}}
\def\fB{{\mathfrak{B}}}
\def\fc{{\mathfrak{c}}}
\def\fC{{\mathfrak{C}}}

\def\fm{{\mathfrak{m}}}

\def\fp{{\mathfrak{p}}}

\def\fq{{\mathfrak{q}}}

\def\nn{\relax\ifmmode{\mathbb N_{0}}\else$\mathbb N_{0}$\fi}
\def\lra{\longrightarrow}

\begin{document}

\title[TIGHT CLOSURE TEST EXPONENTS]{TIGHT CLOSURE TEST EXPONENTS FOR CERTAIN
PARAMETER IDEALS}
\author{RODNEY Y. SHARP}
\address{Department of Pure Mathematics,
University of Sheffield, Hicks Building, Sheffield S3 7RH, United Kingdom\\
{\it Fax number}: 0044-114-222-3769}
\email{R.Y.Sharp@sheffield.ac.uk}

\thanks{The author was partially supported by the
Engineering and Physical Sciences Research Council of the United
Kingdom.}

\subjclass[2000]{Primary 13A35, 13A15, 13D45, 13E05, 13E10, 13H10,
16S36; Secondary 13C15}

\date{\today}

\keywords{Commutative Noetherian ring, prime characteristic,
Frobenius homomorphism, Artinian module, skew polynomial ring,
tight closure, weak test element, test exponent, parameter ideal.}

\begin{abstract}
This paper is concerned with the tight closure of an ideal $\fa$
in a commutative Noetherian ring $R$ of prime characteristic $p$.
The formal definition requires, on the face of things, an infinite
number of checks to determine whether or not an element of $R$
belongs to the tight closure of $\fa$. The situation in this
respect is much improved by Hochster's and Huneke's test elements
for tight closure, which exist when $R$ is a reduced algebra of
finite type over an excellent local ring of characteristic $p$.

More recently, Hochster and Huneke have introduced the concept of
test exponent for tight closure: existence of these test exponents
would mean that one would have to perform just one single check to
determine whether or not an element of $R$ belongs to the tight
closure of $\fa$. However, to quote Hochster and Huneke, `it is
not at all clear whether to expect test exponents to exist;
roughly speaking, test exponents exist if and only if tight
closure commutes with localization'.

The main purpose of this paper is to provide a short direct proof that test
exponents exist for parameter ideals in a reduced excellent
equidimensional local ring of characteristic $p$.
\end{abstract}

\maketitle

\setcounter{section}{-1}
\section{\sc Introduction}
\label{in}

Throughout the paper, $R$ will denote a commutative
Noetherian ring of prime characteristic $p$.
We shall always
denote by $f:R\lra R$ the Frobenius homomorphism, for which $f(r)
= r^p$ for all $r \in R$. Let $\fa$ be an ideal of $R$. The {\em
$n$-th Frobenius power\/} $\fa^{[p^n]}$ of $\fa$ is the ideal of
$R$ generated by all $p^n$-th powers of elements of $\fa$.

We use $R^{\circ}$ to denote the complement in $R$ of the union of
the minimal prime ideals of $R$. An element $r \in R$ belongs to
the {\em tight closure $\fa^*$ of $\fa$\/} if and only if there
exists $c \in R^{\circ}$ such that $cr^{p^n} \in \fa^{[p^n]}$ for
all $n \gg 0$. We say that $\fa$ is {\em tightly closed\/}
precisely when $\fa^* = \fa$. The theory of tight closure was
invented by M. Hochster and C. Huneke \cite{HocHun90}, and many
applications have been found for the theory: see \cite{Hunek96}.
For the definition of the tight closure $N^*_M$ of a submodule $N$
in an ambient $R$-module $M$ (and explanation of the notations
$N^{[p^n]}_M$ and $m^{p^n}$ for $m \in M$ and a non-negative
integer $n$), the reader is referred to \cite[(8.1), (8.2) and
(8.3)]{HocHun90}.

A {\em $p^{w_0}$-weak test element\/} for $R$ (where $w_0$ is a
non-negative integer) is an element $c' \in R^{\circ}$ such that,
for every finitely generated $R$-module $M$ and every submodule
$N$ of $M$, and for $m \in M$, it is the case that $m \in N^*_M$
if and only if $c'm^{p^n} \in N^{[p^n]}_M$ for all $n \geq w_0$. A
$p^0$-weak test element is called a {\em test element\/}. A {\em
locally stable $p^{w_0}$-weak test element\/} (respectively {\em
completely stable $p^{w_0}$-weak test element\/}) for $R$ is an
element $c' \in R$ such that, for every prime ideal $\fp$ of $R$,
the natural image $c'/1$ of $c'$ in the localization $R_{\fp}$ is
a $p^{w_0}$-weak test element for $R_{\fp}$ (respectively, for the
completion $\widehat{R_{\fp}}$ of $R_{\fp}$). When $w_0 = 0$, we
omit the adjective `$p^{w_0}$-weak'. A locally stable
$p^{w_0}$-weak test element for $R$ is a $p^{w_0}$-weak test
element for $R$, and a completely stable $p^{w_0}$-weak test
element for $R$ is a locally stable $p^{w_0}$-weak test element
for $R$: see \cite[Proposition (8.13)]{HocHun90}.

It is a result of Hochster and Huneke \cite[Theorem
(6.1)(b)]{HocHun94} that an algebra of finite type over an
excellent local ring of characteristic $p$ has a completely stable
$p^{w_0}$-weak test element for some $w_0$; furthermore, such an
algebra which is reduced actually has a completely stable test
element.

This paper is concerned with the concept of {\em test exponent\/}
in tight closure theory introduced by Hochster and Huneke in
\cite[Definition 2.2]{HocHun00}. Let $c$ be a test element for a
reduced commutative Noetherian ring $R$ of characteristic $p$,
and let $\fa$ be an ideal of $R$. A test exponent for $c$, $\fa$
is a power $q = p^{e_0}$ (where $e_0 \in \nn$, the set of
non-negative integers) such that if, for an $r \in R$, we have
$cr^{p^e} \in \fa^{[p^e]}$ for {\em one single\/} $e \geq e_0$,
then $r \in \fa^*$ (so that $cr^{p^n} \in \fa^{[p^n]}$ for all $n
\in \nn$). In \cite{HocHun00}, it is shown that this concept has
strong connections with the major open problem about whether tight
closure commutes with localization; indeed, to quote Hochster and
Huneke, `roughly speaking, test exponents exist
if and only if tight closure commutes with localization'.
Although the question of whether tight closure commutes with
localization is open in general, it is known that it does in many
particular cases (see \cite{AHH}); consequently, the results of
Hochster and Huneke in \cite{HocHun00} imply (via a rather
circuitous route) that test exponents must exist rather often.

In \cite[Discussion 5.1]{HocHun00}, Hochster and Huneke state that
`it would be of considerable interest to solve the problem of determining
test exponents effectively even for parameter ideals'.
The main purpose of this paper is to provide a short direct proof
that, for a test element $c$ for a reduced excellent
equidimensional local ring $(R,\fm)$,
there exists $e_0 \in \nn$ such that $p^{e_0}$ is a
test exponent for $c$, $\fa$ for every parameter ideal $\fa$ of
$R$. (For such $(R,\fm)$, a {\em parameter ideal\/} of
$R$ is simply one that can be generated by a subset of a system of
parameters for $R$.) The fact that $p^{e_0}$ is a test exponent
for $c$, $\fa$ for {\em every\/} parameter ideal $\fa$ of $R$ is
relevant to \cite[Discussion 5.3]{HocHun00}, where Hochster and
Huneke raise the question as to whether there might conceivably
exist (when $R$ (not necessarily local) and $c$ satisfy certain
conditions) a `uniform test exponent' for $c$, that is, a power of
$p$ that is a test exponent for $c$, $\fb$ for {\em all\/} ideals
$\fb$ of $R$ simultaneously.

\section{\sc Left modules over the skew polynomial ring $R[x,f]$}
\label{sp}

\begin{ntn}
\label{nt.1}

We shall work with the
 skew polynomial ring $R[x,f]$ associated to $R$ and $f$
in the indeterminate $x$ over $R$. Recall that $R[x,f]$ is, as a
left $R$-module, freely generated by $(x^i)_{i \in \nn}$,
 and so consists
 of all polynomials $\sum_{i = 0}^n r_i x^i$, where  $n \in \nn$
 and  $r_0,\ldots,r_n \in R$; however, its multiplication is subject to the
 rule
 $$
  xr = f(r)x = r^px \quad \mbox{~for all~} r \in R\/.
 $$
Note that the decomposition $R[x,f] = \bigoplus_{n \in \nn}Rx^n$
provides $R$ with a structure as a positively-graded ring.

We use $\N$ to denote the set of positive integers.
\end{ntn}

The first lemma below enables one to see quickly that, in certain
circumstances, an $R$-module $M$ has a structure as left
$R[x,f]$-module extending its $R$-module structure.

\begin{lem}
\label{nt.3} {\rm (See \cite[Lemma 1.3]{KS}, for example.)} Let
$G$ be an $R$-module and let $\xi : G \lra G$ be a
$\Z$-endomorphism of $G$ such that $\xi(rg) = r^p\xi (g)$ for all
$r \in R$ and $g \in G$. Then the $R$-module structure on $G$ can
be extended to a structure of left $R[x,f]$-module in such a way
that $xg = \xi (g)$ for all $g \in G$.
\end{lem}

\begin{defs}
\label{nt.2a} Let $H$ be a left $R[x,f]$-module. The $R[x,f]$-submodule
$$\Gamma_x(H) := \left\{ h \in H : x^jh = 0 \mbox{~for some~} j
\in \N \right\}$$ of $H$ is called the\/
{\em $x$-torsion submodule} of $H$. We say that $H$ is {\em $x$-torsion\/}
precisely when $H = \Gamma_x(H)$, and that $H$ is {\em $x$-torsion-free\/}
precisely when $\Gamma_x(H) = 0$.

It is easy to check that, in general,
the left $R[x,f]$-module $H/\Gamma_x(H)$ is $x$-torsion-free.
\end{defs}

\begin{defs}
\label{hs.1d}
Let $H$ be a left $R[x,f]$-module.
The {\em annihilator of $H$\/}
will be denoted by $\ann_{R[x,f]}H$ or $\ann_{R[x,f]}(H)$. Thus
$
\ann_{R[x,f]}(H) = \{ \theta \in R[x,f] : \theta h = 0 \mbox{~for all~} h \in H\},
$
and this is a (two-sided) ideal of $R[x,f]$.

For a two-sided ideal $\fB$ of $R[x,f]$, we shall use $\ann_H\fB$ or
$\ann_H(\fB)$ to
denote the {\em annihilator of $\fB$ in $H$}. Thus
$
\ann_H\fB = \ann_H(\fB) = \{ h \in
H : \theta h = 0 \mbox{~for all~}\theta \in \fB\},
$
and this is an  $R[x,f]$-submodule of $H$.
\end{defs}

\begin{rmk}
\label{hs.1e}
It is easy to see that a
subset $\fB$ of $R[x,f]$ is a graded left ideal if and only if
there is a family $(\fb_n)_{n\in\nn}$ of ideals
of $R$ such that $\fb_n \subseteq f^{-1}(\fb_{n+1})$ for all $n \in \nn$, and
$\fB = \bigoplus _{n\in\nn} \fb_n x^n$.  Similarly, a
subset $\fC$ of $R[x,f]$ is a graded two-sided ideal if and only if
there is a family $(\fc_n)_{n\in\nn}$ of ideals
of $R$ such that $\fc_n \subseteq \fc_{n+1}$ for all $n \in \nn$ (so that
the sequence $(\fc_n)_{n\in\nn}$ is eventually stationary), and
$\fC = \bigoplus _{n\in\nn} \fc_n x^n$.
\end{rmk}

\begin{lem}
\label{tc.3} Let $G$ be an $x$-torsion-free left $R[x,f]$-module.
Suppose that, for some $w_0 \in \nn$ and some ideal $\fc$ of $R$,
the graded two-sided ideal $\bigoplus_{n \geq w_0} \fc x^n$ of
$R[x,f]$ annihilates $G$. Then $G$ is annihilated by $\bigoplus_{n
\geq 0} \left(\sqrt{\fc}\right) x^n$.
\end{lem}

\begin{proof} Let $a \in R$ be such that $a^t \in \fc$ for some $t \in \N$. We
show that, for $g \in G$, we must have $ax^ng = 0$ for each $n \in
\nn$. Choose $m \in \N$ such that $p^m \geq t$ and $m \geq w_0$;
then $x^max^ng = a^{p^m}x^{m+n}g = 0$ because $a^{p^m} \in \fc$
and $m+n \geq w_0$; since $G$ is $x$-torsion-free, it follows that
$ax^ng = 0$.
\end{proof}

In this paper,
substantial use will be made of the following extension, due to G. Lyubeznik, of
a result of R. Hartshorne and
R. Speiser. It shows that, when $R$ is local, an
$x$-torsion left $R[x,f]$-module which is Artinian
(that is, `cofinite' in the terminology of Hartshorne and Speiser)
as an $R$-module exhibits a certain uniformity of behaviour.

\begin{thm} [G. Lyubeznik {\cite[Proposition 4.4]{Lyube97}}]
\label{hs.4}  {\rm (Compare Hartshorne--Speiser \cite[Proposition
1.11]{HarSpe77}.)} Suppose that $(R,\fm)$ is local, and let $H$ be
a left $R[x,f]$-module which is Artinian as an $R$-module. Then
there exists $e \in \nn$ such that $x^e\Gamma_x(H) = 0$.
\end{thm}

Hartshorne and Speiser first proved this result in the
case where $R$ is local and contains its residue field which is perfect.
Lyubeznik applied his theory of $F$-modules to obtain the result
without restriction on the local ring $R$ of characteristic $p$.

The corollary below extends the Hartshorne--Speiser--Lyubeznik Theorem to
non-local situations.

\begin{cor}
\label{hs.5a} Let
$H$ be a left $R[x,f]$-module which is Artinian as an $R$-module. Then
there exists $e \in \nn$ such that $x^e\Gamma_x(H) = 0$.
\end{cor}

\begin{proof} Suppose that $H \neq 0$.
For each ideal $\fa$ of $R$, let $\Gamma_{\fa}(H) := \bigcup_{n\in\N}(0:_H\fa^n)$.
The ideas of \cite[Exercises 8.48 and 8.49]{SCA} can be used to show that
there are only finitely many maximal ideals $\fm$ of $R$ such that
$\Gamma_{\fm}(H) \neq 0$, and that if we denote the distinct such maximal ideals
by $\fm_1, \ldots \fm_t$, then $H$ decomposes as a direct sum of $R$-submodules
$$
H = \Gamma_{\fm_1}(H) \oplus \cdots \oplus \Gamma_{\fm_t}(H).
$$
Let $i \in \{1, \ldots, t\}$. In fact, $\Gamma_{\fm_i}(H)$ is an $R[x,f]$-submodule
of $H$, since if $h \in \Gamma_{\fm_i}(H)$ and $\fm_i^kh = 0$ for a $k \in \N$,
then $(\fm_i^k)^{[p]}xh = 0$. In addition, for $s \in R \setminus \fm_i$, we have
$\fm_i^k + Rs = R$, so that there exists $s' \in R$ such that $s'sh = h$. It follows
that multiplication by $s$ provides an $R$-automorphism of $\Gamma_{\fm_i}(H)$,
so that the latter left $R[x,f]$-module has a natural structure as an
$R_{\fm_i}$-module in which $(r/s)h$, for $r \in R$ and $s$ as above, is
equal to the unique element $h' \in H$ for which $sh' = rh$. It can
easily be checked that this structure is such that
$
x(r/s)h = (r^p/s^p)xh,
$
and so it follows from Lemma \ref{nt.3} that this $R_{\fm_i}$-module structure on
$\Gamma_{\fm_i}(H)$ can be extended to a structure as left $R_{\fm_i}[x,f]$-module
which is compatible with its structure as a left $R[x,f]$-module. We can now
use the Hartshorne--Speiser--Lyubeznik Theorem \ref{hs.4} to deduce that
there exists $e_i \in \N$ such that
$x^{e_i}\Gamma_x\left(\Gamma_{\fm_i}(H)\right) = 0$.

Since
$
\Gamma_{x}(H)
= \Gamma_{x}\left(\Gamma_{\fm_1}(H)\right) \oplus \cdots \oplus
\Gamma_{x}\left(\Gamma_{\fm_t}(H)\right),
$
the integer $e := \max \{e_1, \ldots,e_t\}$ has the property that
$x^e\Gamma_x(H) = 0$.
\end{proof}

\begin{defi}
\label{hslno}Let $H$ be
a left $R[x,f]$-module which is Artinian as an $R$-module.  By the
Corollary \ref{hs.5a} to the
Hartshorne--Speiser--Lyubeznik Theorem \ref{hs.4},
there exists $e \in \nn$ such that $x^e\Gamma_x(H) = 0$: we call
the smallest such $e$ the {\em Hartshorne--Speiser--Lyubeznik number\/},
or {\em HSL-number\/} for short, of $H$.
\end{defi}

\begin{lem}
\label{hs.3} Let $H$ be
a left $R[x,f]$-module which is Artinian as an $R$-module; let $m_0$ be the
HSL-number of $H$ (see\/ {\rm \ref{hslno}}). Let $\fc$ be an ideal of $R$ and
let $t_0 \in \nn$.

Then $\ann_H\!\left(\bigoplus_{n \geq m_0 + t_0} \fc^{[p^{m_0}]}x^n\right)$ is
an $R[x,f]$-submodule of $H$ that contains $\Gamma_x(H)$; furthermore,
\begin{align*}
\left(\ann_H\!\left({\textstyle \bigoplus_{n \geq m_0 + t_0}
\fc^{[p^{m_0}]}x^n}\right)\right)\! \Big/\Gamma_x(H) &=
\ann_{H/\Gamma_x(H)}\!\left({\textstyle \bigoplus_{n \geq t_0} \fc
x^n}\right)\\ &= \ann_{H/\Gamma_x(H)}\!\left({\textstyle
\bigoplus_{n \geq 0}}
\left(\sqrt{\fc}\right) x^n\right) \\
&= \left(\ann_H\!\left({\textstyle \bigoplus_{n \geq m_0}}
\left(\sqrt{\fc}
\right)^{[p^{m_0}]}x^n\right)\right)\!\Big/\Gamma_x(H).
\end{align*}
\end{lem}

\begin{proof} Since
$\bigoplus_{n \geq m_0 + t_0} \fc^{[p^{m_0}]}x^n \subseteq
\bigoplus_{n \geq m_0} Rx^n$, it is immediate that
$$
\Gamma_x(H) = \ann_H\!\left({\textstyle \bigoplus_{n \geq m_0}}
Rx^n\right) \subseteq \ann_H\!\left({\textstyle \bigoplus_{n \geq
m_0 + t_0}} \fc^{[p^{m_0}]}x^n\right).
$$
Now let $h \in H$, $n \in \nn$ and $c \in \fc$. Then $cx^n(h + \Gamma_x(H)) = 0$
in $H/\Gamma_x(H)$ if and only
if $cx^nh \in \Gamma_x(H)$, and, by definition of the HSL-number of $H$, this is
the case if and only if $x^{m_0}cx^nh = 0$ in $H$, that is, if and only if
$c^{p^{m_0}}x^{m_0+n}h = 0$. It is now easy to prove all the claims by means
of this observation and Lemma \ref{tc.3}.
\end{proof}

The next theorem is the key result of this paper.

\begin{thm}
\label{hs.5} Let $G$ be
an $x$-torsion-free left $R[x,f]$-module which is Artinian as an $R$-module;
let $\fc$ be an ideal of $R$. Let $N$ be the $R$-submodule $(0:_G\fc)$ of $G$; for
each $i \in \nn$, set
$$
N_i := \left\{ g \in G: x^ig \in N\right\} = \left\{ g \in G: c x^ig = 0
\mbox{~for all~} c \in \fc \right\}
= \left\{ g \in G: \fc x^ig = 0\right\}.
$$

\begin{enumerate}
\item Each $N_i~(i \in \nn)$ is an $R$-submodule of $G$.
\item We have $N_i \supseteq N_{i+1}$ for all $i \in \nn$.
\item If $N_i = N_{i+1}$ for some $i \in \nn$, then $N_{i+1} = N_{i+2}$.
\end{enumerate}
Since $G$ is Artinian as an $R$-module, it follows from\/ {\rm
(i), (ii)} and\/ {\rm (iii)} above that there exists (a uniquely
determined) $v_0 \in \nn$ such that
$$
N = N_0 \supset N_1 \supset \cdots \supset N_{v_0} =  N_{v_0+1} = \cdots =
N_{v_0+j} = \cdots
$$
(where `$\supset$' is reserved to denote strict containment). Then $N_{v_0}$
is the largest $R[x,f]$-submodule of $G$ that is contained in $N$.

Henceforth, we shall refer to
the integer $v_0$ as\/ {\em the $\fc$-stability index of $G$}.
Note that it has the following property: for $g \in G$, if
$\fc x^{n_1}g = 0$ for one single integer $n_1 \geq v_0$,
then $\fc x^ng = 0$ for all $n \in \nn$.
\end{thm}

\begin{proof} (i) It is clear that $N_i$ is an Abelian subgroup of $G$. Let
$g \in N_i$ and let $r \in R$. Thus $cx^ig = 0$ for all $c \in \fc$. Hence
$cx^i(rg) = cr^{p^i}x^ig = r^{p^i}cx^ig = 0$ for all $c \in \fc$, so that $rg \in
N_i$.

(ii) Let $g \in N_{i+1}$, and let $c \in \fc$. Therefore $cx^{i+1}g = 0$. Hence
$xcx^ig = c^px^{i+1}g = 0$; since $G$ is $x$-torsion-free, it follows that
$cx^ig = 0$. Therefore $g \in N_i$.

(iii) Assume that $N_i = N_{i+1}$. By part (ii), we have $N_{i+1}
\supseteq N_{i+2}$. Let $g \in N_{i+1}$ and let $c \in \fc$; thus
$cx^{i+1}g = 0$. Therefore $cx^i(xg) = 0$, so that (as this is
true for all $c \in \fc$) we must have $xg \in N_i = N_{i+1}$.
Therefore $cx^{i+1}(xg) = 0$, that is, $cx^{i+2}g = 0$, for all $c
\in \fc$. Hence $g \in N_{i+2}$.

The only remaining claim that still requires proof is the one that
$N_{v_0}$ is the largest $R[x,f]$-submodule of $G$ that is
contained in $N$. To see this, first note that, if $g \in
N_{v_0}$, then $g \in N_{v_0+1}$, so that $cx^{v_0}(xg) =
cx^{v_0+1}g = 0$ for all $c \in \fc$ and $xg \in N_{v_0}$. This
shows that $N_{v_0}$ is an $R[x,f]$-submodule of $G$; it is
contained in $N = N_0$ by part (ii) above. On the other hand, if
$L$ is any $R[x,f]$-submodule of $G$ that is contained in $N$, and
$g \in L$, then we must have $x^{v_0}g \in L \subseteq N$, so that
$g \in N_{v_0}$; hence $L \subseteq N_{v_0}$.
\end{proof}

Next, we use the Corollary \ref{hs.5a} to the
Hartshorne--Speiser--Lyubeznik Theorem \ref{hs.4} to produce a
consequence of Theorem \ref{hs.5} that applies to a left
$R[x,f]$-module that is Artinian as an $R$-module but which is not
necessarily $x$-torsion-free.

\begin{cor}
\label{hs.7} Let $H$ be a left $R[x,f]$-module which is Artinian
as an $R$-module, and let $m_0$ be its HSL-number (see\/ {\rm
\ref{hslno}}). Let $\fc$ be an ideal of $R$, and let $v_0$ be the
$\fc$-stability index (see\/ {\rm \ref{hs.5}})
of the $x$-torsion-free left $R[x,f]$-module
$G := H/\Gamma_x(H)$. Let $h \in H$. Then the following statements
are equivalent:
\begin{enumerate}
\item there exists one single integer $n_1 \geq m_0 + v_0$ such
that $\fc^{[p^{m_0}]}x^{n_1}h = 0$; \item $\fc^{[p^{m_0}]}x^nh =
0$ for all $n \geq m_0$.
\end{enumerate}
\end{cor}

\begin{proof} Suppose that $n_1 \in \nn$ is such that
$n_1 \geq m_0 + v_0$ and $\fc^{[p^{m_0}]}x^{n_1}h = 0$.
Then $x^{m_0}\fc x^{n_1-m_0}h = 0$,
so that $\fc x^{n_1-m_0}h
\subseteq \Gamma_x(H)$ and $\fc x^{n_1-m_0}(h + \Gamma_x(H)) = 0$
in $G$. Since $n_1-m_0 \geq v_0$, the $\fc$-stability index of $G$, it follows
from Theorem \ref{hs.5} that $\fc x^n(h + \Gamma_x(H)) = 0$ for all $n \in \nn$.
Therefore, by Lemma \ref{hs.3},
$$
h + \Gamma_x(H) \in \ann_{H/\Gamma_x(H)}\!\left({\textstyle
\bigoplus_{n \geq 0} \fc x^n}\right) =
\left(\ann_H\!\left({\textstyle\bigoplus_{n \geq m_0}}
\left(\sqrt{\fc}
\right)^{[p^{m_0}]}x^n\right)\right)\!\Big/\Gamma_x(H).
$$
Therefore $h$ is annihilated by $\bigoplus_{n \geq m_0}
\fc^{[p^{m_0}]}x^n$.
\end{proof}

\section{\sc Applications to test exponents for tight closure}
\label{tc}

The  main strategy employed in this paper involves application of the
key result, and its corollary, of \S \ref{sp} to the top local cohomology
module of $(R,\fm)$ in the case when the latter is an equidimensional excellent
local ring (of characteristic $p$).
We therefore review the $R[x,f]$-module structure carried by this
local cohomology module.

\begin{rmd}
\label{lc.1} Suppose that $(R,\fm)$ is a local ring of dimension $d > 0$.
In this reminder, we shall sometimes use
$R'$ to denote $R$ regarded as an $R$-module by means of $f$.

\begin{enumerate}
\item With this notation, $f : R \lra R'$ becomes a homomorphism of $R$-modules, and
so induces an $R$-homomorphism
$
H^d_{\fm}(f) : H^d_{\fm}(R) \lra H^d_{\fm}(R').
$
The Independence Theorem for local cohomology (see \cite[4.2.1]{LC}) applied
to the ring homomorphism $f : R \lra R$ yields an $R$-isomorphism
$
\nu^d_{R}: H^d_{\fm}(R') \stackrel{\cong}{\lra} H^d_{\fm^{[p]}}(R),
$ where
$H^d_{\fm^{[p]}}(R)$ is regarded as an $R$-module via $f$.
Since $\fm$ and $\fm^{[p]}$ have the same radical, $H^d_{\fm}$ and
$H^d_{\fm^{[p]}}$ are the same functor.
Composition yields a $\Z$-endomorphism
$
\xi := \nu^d_{R} \circ H^d_{\fm}(f) : H^d_{\fm}(R) \lra H^d_{\fm}(R)
$
which is such that $\xi(r\gamma) = r^p\xi(\gamma)$ for all $\gamma \in H^d_{\fm}(R)$
and $r \in R$. It therefore follows from Lemma \ref{nt.3} that $H^d_{\fm}(R)$
has a natural structure as left $R[x,f]$-module in which $x\gamma = \xi(\gamma)$
for all $\gamma \in H^d_{\fm}(R)$.

\item It is important to note that this $R[x,f]$-module
structure on $H^d_{\fm}(R)$ does not depend on any choice
of system of parameters for $R$. The reader might like to
consult \cite[2.1]{KS} for amplification of this point.

\item Let $a_1, \ldots, a_d$ be a system of parameters for $R$, and represent
$H^d_{\fm}(R)$ as the $d$-th cohomology module of the \u{C}ech
complex of $R$ with respect to $a_1, \ldots, a_d$, that is, as the
residue class module of $R_{a_1 \ldots a_d}$ modulo the image,
under the \u{C}ech complex `differentiation' map, of
$\bigoplus_{i=1}^dR_{a_1 \ldots a_{i-1}a_{i+1}\ldots a_d}$. See
\cite[\S 5.1]{LC}. We use `$\left[\phantom{=} \right]$' to denote natural
images of elements of $R_{a_1\ldots a_d}$ in this residue class
module.

It is worth noting that, for $r \in R$ and $n \in \nn$, we have
$$
\left[\frac{r}{(a_1\ldots a_d)^n}\right] = 0 \quad \mbox{~in~} H^d_{\fm}(R)
$$
if and only if there exists $k \in \nn$ such that $(a_1\ldots a_d)^kr \in
(a_1^{n+k}, \ldots,a_d^{n+k})R$.

\item The left $R[x,f]$-module structure on $H^d_{\fm}(R)$ is such
that
$$
x\left[\frac{r}{(a_1\ldots a_d)^n}\right] =
\left[\frac{r^p}{(a_1\ldots a_d)^{np}}\right] \quad \mbox{~for all~} r \in R
\mbox{~and~} n \in \nn.
$$
The reader might like to consult \cite[2.3]{KS} for more details.
\end{enumerate}
\end{rmd}

In \cite[Definition 2.2]{HocHun00}, Hochster and Huneke defined
the concept of test exponent in a {\it reduced\/} commutative
Noetherian ring $R$ of prime characteristic. They also defined the
concept for $c$, $N$, $M$, where ($c$ is a test element for $R$
and) $N$ is a submodule of the finitely generated $R$-module $M$.
However, the definition for modules is not pursued in this paper,
which is concerned with test exponents for $c$, $\fa$, where $\fa$
is an ideal of $R$ (and `$R$' is understood to be the ambient
module for $\fa$). On the other hand, there are advantages in
extending the concept to weak test elements in non-reduced rings.

\begin{defi}
\label{tc.1} Let $c$ be a $p^{w_0}$-weak test element (where $w_0
\in \nn$) for the (not necessarily reduced) ring $R$, and let
$\fa$ be an ideal of $R$. We say that $p^{e_0}$ (where $e_0 \in
\nn$) is a {\em test exponent for $c$, $\fa$\/} if, whenever $r
\in R$ is such that $cr^{p^{e}} \in \fa^{[p^{e}]}$ for {\em
one single\/} $e \geq e_0$, then $r \in \fa^*$ (so that
$cr^{p^n} \in \fa^{[p^n]}$ for all $n \geq w_0$).
\end{defi}

Recall that a {\em parameter ideal\/} in a commutative Noetherian
ring is a proper ideal of height $h$ that can be generated by $h$
elements (for some $h \in \nn$). In an equidimensional catenary
local ring, an ideal is a parameter ideal if and only if it can be
generated by a subset of a system of parameters.

\begin{thm}
\label{tc.4} Let $(R,\fm)$ (as in\/ {\rm \ref{nt.1}}) be an
equidimensional excellent local ring of dimension $d > 0$.

By\/ {\rm \ref{lc.1}},
the Artinian $R$-module $H := H^d_{\fm}(R)$ has a natural structure as
a left $R[x,f]$-module; let $m_0$ be its HSL-number (see\/ {\rm
\ref{hslno}}). Let $c \in R^{\circ}$ and let $v_0$ be the
$Rc$-stability index (see\/ {\rm \ref{hs.5}})
of the $x$-torsion-free left $R[x,f]$-module
$G := H/\Gamma_x(H)$.

Then for each parameter ideal $\fa$ of $R$ the following is true: whenever
$r \in R$ is such that $cr^{p^{n_1}} \in \fa^{[p^{n_1}]}$
for one single $n_1 \geq m_0 + v_0$, then $r \in \fa^*$.
\end{thm}

\begin{proof}
We shall first prove the claim when
$\fa$ is an ideal $\fq$ of $R$ generated by a full system of parameters
$a_1, \ldots, a_d$ for $R$.
Use the representation of $H$ as
the $d$-th cohomology module of the \u{C}ech complex of
$R$ with respect to $a_1, \ldots, a_d$ recalled in \ref{lc.1}(iii), and write $a$
for the product $a_1\ldots a_d$.

Set
$
\zeta := \left[r/a\right] \in H.
$
The assumption that $cr^{p^{n_1}} \in \fq^{[p^{n_1}]}$
for an $r \in R$ and an $n_1 \geq m_0 + v_0$ implies (by \ref{lc.1}(iv)) that
$
cx^{n_1}\zeta = \left[cr^{p^{n_1}}/a^{p^{n_1}}\right] = 0.
$
Therefore $Rc^{p^{m_0}}x^{n_1}\zeta = 0$, and
so it follows from Corollary \ref{hs.7} that
$Rc^{p^{m_0}}x^n\zeta = 0$ for all $n \geq m_0$. Hence
$$
\left[\frac{c^{p^{m_0}}r^{p^{n}}}{a^{p^{n}}}\right]
= c^{p^{m_0}}x^{n}\zeta = 0 \quad \mbox{~for all~} n \geq m_0.
$$
It now follows from \ref{lc.1}(iii) that, for all $n \geq m_0$,
there exists $k(n) \in \nn$
such that
$$
c^{p^{m_0}}r^{p^{n}}(a_1\ldots a_d)^{k(n)} \in (a_1^{p^{n}+k(n)},
\ldots, a_d^{p^{n}+k(n)})R.
$$
The next part of the argument is due to K. E. Smith: see the proof of
\cite[Proposition 3.3(i)]{S94}. By repeated use of the colon-capturing
properties of tight closure described in \cite[Theorem 2.9]{S94}, it follows that
$$
c^{p^{m_0}}r^{p^{n}} \in (\fq^{[p^n]})^* \quad \mbox{~for all~} n \geq m_0.
$$
Since $R$ is an excellent local ring, it has a $p^{w_0}$-weak test
element $c'$, for some $w_0 \in \nn$, by Hochster--Huneke
\cite[Theorem 6.1(b)]{HocHun94}. Consequently,
$$
c'(c^{p^{m_0}}r^{p^{n}})^{p^{w_0}} \in (\fq^{[p^n]})^{[p^{w_0}]}
\quad \mbox{~for all~} n \geq m_0,
$$
that is $c'c^{p^{m_0+w_0}}r^{p^{n+w_0}} \in \fq^{[p^{n+w_0}]}$ for all $n \geq m_0$.
Since $c'c^{p^{m_0+w_0}} \in R^{\circ}$, we see that $r \in \fq^*$.

It remains to extend the result to an arbitrary parameter ideal $\fa$ of $R$. Since
$R$ is equidimensional and excellent, there exists a full system of parameters
$u_1, \ldots, u_d$ for $R$ and $i \in \{0,1, \ldots,d\}$ such that
$\fa = (u_1, \ldots, u_i)R$. Suppose that $r \in R$
is such that $cr^{p^{n_1}} \in \fa^{[p^{n_1}]}$
for an $n_1 \geq m_0 + v_0$. Then, for all $t \in \N$,
we have $cr^{p^{n_1}} \in (u_1, \ldots, u_i,u_{i+1}^t, \ldots,
u_d^t)^{[p^{n_1}]}$, so that $r \in (u_1, \ldots, u_i,u_{i+1}^t, \ldots,
u_d^t)^*$ by the first part of this proof. We can now use the
$p^{w_0}$-weak test element
$c'$ to deduce that
$$
c'r^{p^n} \in (u_1^{p^n}, \ldots, u_i^{p^n},u_{i+1}^{p^nt}, \ldots,
u_d^{p^nt})R \quad \mbox{~for all~} n \geq w_0 \mbox{~and all~} t \in \N.
$$
Therefore, by Krull's Intersection Theorem,
$$
c'r^{p^n} \in (u_1^{p^n}, \ldots, u_i^{p^n})R \quad \mbox{~for all~} n \geq w_0,
$$
so that $r \in \fa^*$.
\end{proof}

If, in Theorem \ref{tc.4}, we take the element $c$ to be a weak
test element for $R$, we can immediately deduce the existence of a
test exponent for $c$, $\fa$ for each parameter ideal of $R$; it
should be noted that this test exponent is `partially uniform' in
the sense that the one test exponent for $c$ works for {\em
every\/} parameter ideal of $R$. These results are recorded in
part (i) of Corollary \ref{tc.8} below. I am very grateful to the
referee for her/his suggestion of part (ii) of \ref{tc.8}, which
shows, loosely speaking, that a slightly higher power of $p$ is
not only a test exponent for $c$, $\fa$ for all parameter ideals
$\fa$ of $R$ simultaneously, but also that one need only check
that the `ideal membership test' is satisfied `up to tight
closure'.

\begin{cor}
\label{tc.8} Let $(R,\fm)$ (as in\/ {\rm \ref{nt.1}}) be an
equidimensional excellent local ring of dimension $d > 0$.
Let $c$ be a $p^{w_0}$-weak test element (where $w_0 \in \nn$)
for $R$.

As in Theorem\/ {\rm \ref{tc.4}}, let $m_0$ be the HSL-number (see\/ {\rm
\ref{hslno}}) of $H := H^d_{\fm}(R)$, and let $v_0$ be the
$Rc$-stability index (see\/ {\rm \ref{hs.5}})
of $G := H/\Gamma_x(H)$.

\begin{enumerate}
\item Then $p^{m_0 + v_0}$ is a test exponent for $c$, $\fa$ for all
parameter ideals $\fa$ of $R$ simultaneously.
\item The power $p^{m_0 + v_0 + 1}$ has the following property: for each
parameter ideal $\fa$ of $R$, whenever
$r \in R$ is such that $cr^{p^{n_1}} \in (\fa^{[p^{n_1}]})^*$
for one single $n_1 \geq m_0 + v_0 + 1$, then $r \in \fa^*$.
\end{enumerate}
\end{cor}

\begin{proof} Part (i) is immediate from Theorem \ref{tc.4}, and so we prove
(ii).

We first show that the $Rc^{p^{w_0}+1}$-stability index $v_1$ of
$G$ satisfies $v_1 \leq v_0 + w_0 + 1$. By Theorem \ref{hs.5}, it
is enough, in order to prove this inequality, to show that
$c^{p^{w_0}+1}x^{v_0 + w_0 + 1}g = 0$ (for $g \in G$) implies that
$c^{p^{w_0}+1}x^ng = 0$ for all $n \in \nn$. However
$c^{p^{w_0}+1}x^{v_0 + w_0 + 1}g = 0$ implies that
$c^{p^{w_0+1}}x^{v_0 + w_0 + 1}g = 0$, so that
$$
x^{w_0 + 1}cx^{v_0}g = c^{p^{w_0 + 1}}x^{v_0 + w_0 + 1}g = 0,
$$
and $cx^{v_0}g = 0$ because $G$ is $x$-torsion-free; since $v_0$
is the $Rc$-stability index of $G$, this implies that $cx^ng = 0$
for all $n \in \nn$, so that $c^{p^{w_0}+1}x^ng = 0$ for all $n
\in \nn$. Therefore $v_1 \leq v_0 + w_0 + 1$.

Now suppose that $\fa$ is a parameter ideal of $R$, and $r \in R$
is such that $cr^{p^{n_1}} \in (\fa^{[p^{n_1}]})^*$ for one single
$n_1 \geq m_0 + v_0 + 1$. Since $c$ is a $p^{w_0}$-weak test
element for $R$, we have $c(cr^{p^{n_1}})^{p^{w_0}} \in
(\fa^{[p^{n_1}]})^{[p^{w_0}]}$, that is
$c^{p^{w_0}+1}r^{p^{n_1+w_0}} \in \fa^{[p^{n_1+w_0}]}$. Now $n_1 +
w_0 \geq m_0 + v_0 + w_0 + 1 \geq m_0 + v_1$, and so it follows
from Theorem \ref{tc.4} that $r \in \fa^*$.
\end{proof}

In the final theorem of the paper, we deduce a non-local result
from Corollary \ref{tc.8}. I am again very grateful to the referee
for her/his suggestions that have led to improvements in Theorem
\ref{tc.9}. Note that $R$ is said to be {\em locally
equidimensional\/} precisely when the localization $R_{\fp}$ is
equidimensional for every prime ideal $\fp$ of $R$. The reader is
referred to Aberbach--Hochster--Huneke \cite[p.\ 87]{AHH} for an
explanation of what it means to say that $R$ is {\em of acceptable
type}.

\begin{thm}
\label{tc.9} Let $R$ (as in\/ {\rm \ref{nt.1}}) be a locally
equidimensional ring of acceptable type. Suppose that there exists
a completely stable $p^{w_0}$-weak test element $c$ for $R$ (where
$w_0 \in \nn$). (These conditions would all be satisfied if $R$
was an integral domain and an algebra of finite type over an
excellent local ring of characteristic $p$, by\/ {\rm
\cite[Proposition (5.4)]{AHH}} and\/ {\rm \cite[Theorem
(6.1)(b)]{HocHun94}}.)

Let $\fa$ be a parameter ideal of $R$ of positive height, and let
${\mathcal P} := \{ \fp_1, \ldots, \fp_t\}$
be a finite set of prime ideals of $R$ of positive height such
that
$$ \bigcup_{\fp \in \ass \fa^*} \fp
\subseteq \bigcup_{i=1}^t \fp_i.$$ (For example, ${\mathcal P}$
could be $\ass \fa^*$.)

For each $i = 1, \ldots, t$, let $m_i$ denote the HSL-number
(see\/ {\rm \ref{hslno}}) of the top local cohomology module $H_i
:= H^{\height \fp_i}_{\fp_i R_{\fp_i}}(R_{\fp_i})$ of the local
ring $R_{\fp_i}$, and let $v_i$ denote the
$R_{\fp_i}(c/1)$-stability index (see\/ {\rm \ref{hs.5}}) of
$H_i/\Gamma_x(H_i)$. Then
$$
u_0 := \max\{m_1 + v_1, \ldots, m_t + v_t\}
$$
has the following property: whenever $r \in R$ is such that
$cr^{p^{n_1}} \in (\fa^{[p^{n_1}]})^*$ for one single $n_1 \geq
u_0 + 1$, then $r \in \fa^*$.

Consequently, $p^{u_0 + 1}$ is a test exponent for $c$, $\fa$.
\end{thm}

\begin{proof} Temporarily, let $A$ be a local commutative Noetherian ring of
characteristic $p$ and positive dimension $d$, and suppose that
$H$ is a left $A[x,f]$-module which is Artinian as an $A$-module.
Then $H$ has a natural structure as a module over the completion
$\widehat{A}$ of $A$ (see \cite[8.2.4]{LC}), and it is easy to use
Lemma \ref{nt.3} to see that this $\widehat{A}$-module structure
on $H$ can be extended to a structure as left
$\widehat{A}[x,f]$-module which is compatible with its structure
as a left $A[x,f]$-module. Thus $\Gamma_x(H)$ is the same whether
calculated over $A$ or $\widehat{A}$, and a similar comment
applies to the HSL-number of $H$. Note also that, for $c' \in A$,
the $Ac'$-stability index of $H/\Gamma_x(H)$, as left
$A[x,f]$-module, is the same as its $\widehat{A}c'$-stability index
as left $\widehat{A}[x,f]$-module.

Again by \cite[8.2.4]{LC}, there is an isomorphism of
$\widehat{A}$-modules $H \cong H\otimes_A\widehat{A}$, and so it
follows from \ref{lc.1}(iv) and the comments in the preceding
paragraph that, for each $i \in \{1, \ldots,t\}$,
the HSL-number of the top local cohomology module
$H_i$ of $R_{\fp_i}$ is equal to the corresponding number for
$\widehat{R_{\fp_i}}$, and that the $R_{\fp_i}(c/1)$-stability
index for $H_i/\Gamma_x(H_i)$ is equal to the corresponding index
(for $\widehat{R_{\fp_i}}(c/1)$).

Since $R$ is of acceptable type, it and all its localizations are
universally catenary; therefore, by Ratliff's Theorem (see
\cite[Theorem 31.7]{HM}), all the $R_{\fp_i}~(i = 1, \ldots,t)$
are formally catenary (see \cite[p.\ 252]{HM}). Since $R$ is
locally equidimensional, it follows that all the
$\widehat{R_{\fp_i}}~(i = 1, \ldots,t)$ are equidimensional; they
are also excellent, because every complete Noetherian local ring
is excellent. Thus Corollary \ref{tc.8}(ii) can be applied over
each $\widehat{R_{\fp_i}}~(i = 1, \ldots,t)$.

For each $i = 1, \ldots, t$, let $\phantom{\fa}^{e_i}$ and
$\phantom{\fa}^{c_i}$ stand for extension and
 contraction with respect to the natural ring homomorphism
$R \lra R_{\fp_i}$.

Let $r \in R$ be such that $cr^{p^{n_1}} \in (\fa^{[p^{n_1}]})^*$
for one single $n_1 \geq u_0 + 1$. Choose $i \in \{1, \ldots,
t\}$. Then, in the local ring $R_{\fp_i}$ and its completion, we
have
$$
\frac{c}{1}\left(\frac{r}{1}\right)^{p^{n_1}} =
\frac{cr^{p^{n_1}}}{1} \in ((\fa^{[p^{n_1}]})^*)^{e_i} \subseteq
((\fa^{[p^{n_1}]})^{e_i})^* = ((\fa^{e_i})^{[p^{n_1}]})^*
\subseteq ((\fa^{e_i}\widehat{R_{\fp_i}})^{[p^{n_1}]})^*.
$$
Now $\fa^{e_i}\widehat{R_{\fp_i}}$, if proper, is a parameter
ideal of $\widehat{R_{\fp_i}}$. Since $n_1 \geq m_i + v_i + 1$, it
therefore follows from Corollary \ref{tc.8}(ii) that $r/1 \in
(\fa^{e_i}\widehat{R_{\fp_i}})^*$.

Since $c$ is a completely stable $p^{w_0}$-weak test element for
$R$, it follows that
$$
\frac{c}{1}\left(\frac{r}{1}\right)^{p^n} \in
\left((\fa^{e_i})^{[p^{n}]}\widehat{R_{\fp_i}}\right) \cap
R_{\fp_i} = (\fa^{e_i})^{[p^{n}]} \quad \mbox{~for all~} n \geq
w_0.
$$
Therefore $r/1 \in (\fa^{e_i})^*$. Since $R$ is locally
equidimensional and of acceptable type, we can use \cite[Theorem
(8.3)(a)]{AHH} to see that localization commutes with tight
closure for the pair $\fa \subseteq R$, and so $r/1 \in
(\fa^*)^{e_i}$ and $r \in (\fa^*)^{e_ic_i}$. This is true for all
$i = 1, \ldots, t$.

However, it follows from elementary facts about primary
decomposition that the hypotheses about $\fp_1, \ldots, \fp_t$
ensure that $\fa^* = \bigcap_{i=1}^t (\fa^*)^{e_ic_i}$, and so it
follows that $r \in \fa^*$, as required.

It is then immediate that $p^{u_0 + 1}$ is a test exponent for
$c$, $\fa$.
\end{proof}

\bibliographystyle{amsplain}

\end{document}